\documentclass[preprint,11pt]{elsarticle}
\usepackage{amsmath,amssymb}
\usepackage[all]{xy}
\usepackage{latexsym}
\usepackage{amsthm,a4wide,color}
\usepackage{amsmath,amscd,verbatim}
\usepackage{hyperref}
\input diagxy

\theoremstyle{plain}
  \newtheorem{thm}{Theorem}[section]
  \newtheorem{lem}[thm]{Lemma}

\theoremstyle{definition}
  \newtheorem*{defn}{Definition}

\begin{document}

\begin{frontmatter}
\title{ Four-body Central Configurations with Adjacent Equal Masses}

\author[S]{Yiyang Deng\corref{cor}}
\ead{dengyiyang@126.com}

\author[Y]{Bingyu Li}
\ead{sbtgvpgf@163.com}

\author[S]{Shiqing Zhang}
\ead{zhangshiqing@msn.com}

\cortext[cor]{Corresponding author.}
\address[S]{Department of Mathematics, Sichuan University, Chengdu 610064, China}
\address[Y]{College of Management Science, Chengdu University of Technology, Chengdu 610059, China}

\date{}

\begin{abstract}
 For any convex non-collinear central configuration of the planar  Newtonian 4-body problem with adjacent equal masses $m_1=m_2\neq m_3=m_4$, with equal lengths for the two diagonals, we prove it must possess a symmetry and must be an isosceles trapezoid; furthermore, which is also an isosceles trapezoid when the length between $m_1$ and $ m_4$ equals the length between $m_2$ and $ m_3$.
\end{abstract}
\begin{keyword} Four-body problem, Convex central configurations, Isosceles trapezoid, Triangle areas \end{keyword}

\end{frontmatter}

2000 Mathematical Subject Classification 34C15, 34C25, 70F07
\section{Introduction}

It is well known\cite{DS80}\cite{W41} that the central configurations for Newtonian n-body problems play an important role in the celestial mechanics.
Here, we only consider the planar Newtonian 4-body problem.

The planar Newtonian 4-body problem is related with the motion of 4 point particles with positive masses $m_i\in \mathbb{R}$ and position vectors $q_i\in \mathbb{R}^2$ for $i=1,\cdots,4$, moving according to Newton's second law and the universal gravitational law:
\begin{equation} m_i\ddot {q_i}=\frac{\partial U(q)}{\partial q_i },\ \ \ i=1,\ldots,4,   \end{equation}
where \begin{equation} U(q)=G \sum_{i<j}^{4} \frac{m_{i}m_{j}}{r_{ij}} \label{First}   \end{equation}
is the Newtonian potential for the 4-body and $r_{ij}=\|q_i-q_j\|$, in the following, we let $G=1$.

Let $q=(q_1,\cdots,q_4)\in (\mathbb{R}^2)^4$ and  $M$ be the diagonal mass matrix ${\rm diag}(m_1,m_1,\cdots,m_4,m_4)$,
then the system $(1)$ can be rewritten as the following:
\begin{equation}\ddot{q}=M^{-1}\frac{\partial U(q)}{\partial (q)}. \end{equation}

To study this problem, without lose of generality, we assume the center of mass is fixed at the origin and consider the space
$$\Omega=\left\{{ q=(q_1,q_2,q_3,q_4)\in (\mathbb{R}^2)^4 \mid \sum_{i=1}^{4} m_{i}q_{i}=0 }\right\}.$$
Let $\Delta=\bigcup_{i\neq j}\{q \mid q_i = q_j \}$ be the \textbf{collision set}.
The set $\Omega \setminus \Delta$ is called the \textbf{configuration space}.

Here we recall the definition of the central configuration:
\begin{defn}{\cite{DS80}\cite{W41}}
A configuration $q\in \Omega \setminus \Delta$ is called a central configuration if there is some constant $\lambda$ such that \begin{equation} M^{-1}\frac{\partial U}{\partial q}=\lambda q. \end{equation}
\end{defn}
This equation is invariant under rotation, dilatation and reflection on the plane. Two cental configurations are considered equivalent if they are related by those symmetry operations.

In 1995 and 1996, Albouy \cite{A95}\cite{A96} proved that there are exactly four equivalent classes for the central configurations of the planar Newtonian 4-body problem with positive equal masses.
In 2002, Long and Sun \cite{L02} showed that any convex non-collinear central configurations of the planar 4-body problem with equal opposite masses $\beta>\alpha>0$, such that the diagonal corresponding to the mass $\alpha$ is not shorter than that corresponding to the mass $\beta$, must possess a symmetry and must be a kite. Furthermore, it must be a rhombus.
In 2003, Albouy \cite{A03} showed that any convex  non-collinear  central configurations of the planar 4-body problem can not be a kite when two pairs of opposite masses are not equal.
In 2007, Perez-Chavela and Santoprete \cite{PC07} generalized the result of Long and Sun\cite{L02} and obtained the symmetry of central configurations with equal masses located at opposite vertices of a quadrilateral, but they assume that the two equal masses are not the smallest in all masses.
In 2008, Abouy, Fu and Sun \cite{A08} proved that, in the planar 4-body problem, a convex central configuration is symmetric with respect to one diagonal if and only if the masses of the two particles on the other diagonal are equal. They also showed that the less massive one is closer to the former diagonal. In this paper, they raised a question:
Does the equality of two pairs of adjacent masses implies the configuration is an isosceles trapezoid for co-planar 4-body convex central configurations?

 \setlength{\unitlength}{5cm}
\begin{center}
\begin{picture}(1.25,1)
\put(0,0){\line(1,0){1.2}}
\put(0,0){\line(1,2){0.4}}
\put(0.4,0.8){\line(1,0){0.6}}
\put(1,0.8){\line(1,-4){0.2}}
\put(0,0){\line(5,4){1}}
\put(0.4,0.8){\line(1,-1){0.8}}
\put(-0.1,0){${q_1}$}
\put(1.2,0){${q_2}$}
\put(1,0.8){${q_3}$}
\put(0.3,0.8){${q_4}$}
\end{picture}

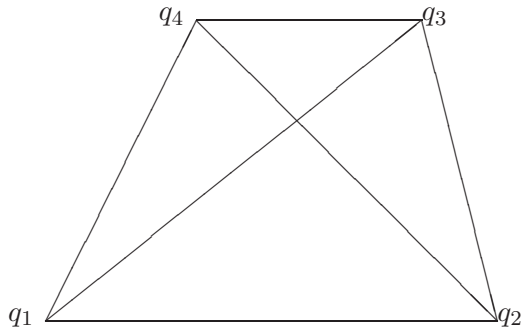
{\bf{Figure.1}}
\end{center}

To solve this problem, in 2012, Cors and Roberts \cite{CR12} used mutual distances as coordinates to study the four-body co-circular central configurations.
They had proved that the set of positions that yield co-circular central configurations with positive masses is a two-dimensional surface, the graph of a differentiable function over two of the exterior side-lengths.
The boundary of this surface correspond to three important symmetric cases: a kite, an isosceles trapezoid and a degenerate case where three bodies lie at the vertices of an equilateral triangle and the fourth body of the quadrilateral has zero mass; furthermore, they got a stronger result, they only assume that one pair of adjacent masses are equal, then the configuration is an isosceles trapezoid for 4-body convex co-circular central configurations.
In 2014, Corbera and Llibre \cite{Cl14} showed that there is a unique convex planar central configuration which has two pairs of equal masses located at the adjacent vertices of the configuration, and which is an isosceles trapezoid when one pair of adjacent masses is sufficiently small.

In this paper, we study the isosceles trapezoid central configurations for the four-body problems. Our main
results are:

\begin{thm}
Let $q=(q_1,q_2,q_3,q_4)\in \Omega$ be a convex non-collinear central configuration with masses $(\beta,\beta,\alpha,\alpha)$, $\beta>\alpha>0$. Suppose that the equal masses are at adjacent vertices. If $$r_{13}=r_{24},$$ then the configuration $q$ must possess a symmetry, and then forms an isosceles trapezoid.\end{thm}

\begin{thm}
Let $q=(q_1,q_2,q_3,q_4)\in \Omega$ be a convex non-collinear central configuration with masses $(\beta,\beta,\alpha,\alpha)$, $\beta>\alpha>0$. Suppose that the equal masses are at adjacent vertices. If $$r_{14}=r_{23},$$ then the configuration $q$ must possess a symmetry, and then forms an isosceles trapezoid.\end{thm}

The arrangement of this paper is as follows.
In section 2, we establish our equations for 4-body central configurations by using areas; In section 3, we will prove Theorem 1.2; In section 4, we will give the proof of Theorem 1.3.

\section{The central configurations equations by areas of triangles}
Firstly, we observe that if  $q=(q_1,q_2,q_3,q_4)\in \Omega$ is a central configuration with parameter $\lambda$ and positive masses $(m_1,m_2,m_3,m_4)$, then every $\zeta^{-\frac{1}{3}}(q_1,q_2,q_3,q_4)\in \Omega$ is in the same class of central configurations with masses $\zeta^{-1}(m_1,m_2,m_3,m_4)$ and the same value of $\lambda$. So, without loss of generality, we suppose $\beta=1$, and we consider the planar 4-body problem with masses $$m_1=m_2=1,~~~~~m_3=m_4=\alpha.$$

In this paper, we use Dziobeck coordinates, which will be described below. Let $$a=r_{12}^2,~b=r_{13}^2,~c=r_{14}^2,~d=r_{23}^2,~e=r_{24}^2,~f=r_{34}^2.$$

For $1\leqq i \leqq4$, let $|\triangle_i|$ be the area of the sub-triangle formed by the remaining three vertices of the configuration $q$ when deleting the point $q_i$. Then we define the oriented areas of these sub-triangles of the convex non-collinear configuration $q$ by
\begin{equation}
\triangle_1=-|\triangle_1|, \ \triangle_2=|\triangle_2|, \ \triangle_3=-|\triangle_3|,\ \triangle_4=|\triangle_4|.
\end{equation}
The above $\triangle_i$ satisfy the following equality:
\begin{equation}\triangle_1+\triangle_2+\triangle_3+\triangle_4=0.\end{equation}

It is well known\cite{D19}\cite{DS02} that the Cayley-Menger determinant
   $$ S=\begin{vmatrix}
       0 & 1 & 1 & 1 & 1 \\
       1 & 0 & a & b & c \\
       1 & a & 0 & d & e \\
       1 & b & d & 0 & f \\
       1 & c & e & f & 0 \\
     \end{vmatrix}$$
satisfies $S=0$. In 1900, Dziobek \cite{D19}(also refer to\cite{CR12},\cite{M01}) proved that
\begin{equation}\frac{\partial S}{\partial r^2_{ij}}=-32\triangle_i \triangle_j,\ \ \ \rm{for~ all} \ i\neq j.\end{equation}

Let $\varphi(s)=s^{-\frac{1}{2}}$ for $s>0$. Then the potential function and the momentum of inertia are given by
\begin{equation}U(q)=\sum_{1\leqq i < j\leqq 4 }{m_i}{m_j}{\varphi(r^2_{ij})} \end{equation}   and
\begin{equation}I(q)=\frac{1}{m'}\sum_{1\leqq i < j\leqq 4 }{m_i}{m_j}{r^2_{ij}}\end{equation}
respectively, where $m'=\sum_{i=1}^{4}m_i$.

Using Lagrangian Multiplier Method, Dziobek gave an equivalent characterization of central configurations, they are extremal of $$U+\lambda S-\mu m'(I-I_0)$$ as a function of $\lambda ,\ \mu m',\ r_{12},\cdots,r_{34}$, where $\lambda $ and $\mu m'$ are Lagrange multipliers and $I_0$ is a fixed moment of inertia. Thus, for any $i,j$ with $1\leqq i<j\leqq 4$, the central configuration satisfies
 \begin{equation}\frac{\partial U}{\partial r^2_{ij}}=-\lambda \frac{\partial S}{\partial r^2_{ij}}+\mu m'\frac{\partial I}{\partial r^2_{ij}}.\end{equation}

 By $(8)$ and $(9)$, we have $$\frac{\partial U}{\partial r^2_{ij}}={m_i}{m_j}{\varphi'(r^2_{ij})},$$
where $\varphi'(s)$ denotes the derivative of function $\varphi(s)$ with respect to $s$, and
$$m'\frac{\partial I}{\partial r^2_{ij}}=m_i m_j.$$
So, the equation $(10)$ becomes
\begin{equation} {m_i}{m_j}{\varphi'(r^2_{ij})}=32\lambda \triangle_i \triangle_j +\mu {m_i}{m_j}. \end{equation}

Using our assumption on masses, the equations of the central configurations become
\begin{equation} \varphi'(r^2_{12})= \nu \triangle_1 \triangle_2 + \mu, \end{equation}
\begin{equation} \varphi'(r^2_{13})= \frac{\nu}{\alpha} \triangle_1 \triangle_3 + \mu, \end{equation}
\begin{equation} \varphi'(r^2_{14})= \frac{\nu}{\alpha} \triangle_1 \triangle_4 + \mu, \end{equation}
\begin{equation} \varphi'(r^2_{23})= \frac{\nu}{\alpha} \triangle_2 \triangle_3 + \mu, \end{equation}
\begin{equation} \varphi'(r^2_{24})= \frac{\nu}{\alpha} \triangle_2 \triangle_4 + \mu, \end{equation}
\begin{equation} \varphi'(r^2_{34})= \frac{\nu}{\alpha^2}\triangle_3 \triangle_4 + \mu, \end{equation}
where $\nu=32\lambda$.

According to Albouy\cite{A95}\cite{A96}, the geometrical relations between $r^2_{ij}$ and $\triangle_i$ are in the following
\begin{equation} t_l=\sum_{i=1}^4 \triangle_i r^2_{il},\ \ t_1=t_2=t_3=t_4. \end{equation}

Using the above implicit relations, Long and Sun, Perez-Chavela and Santoprete got the following Lemma:
\begin{lem}{\cite{L02}\cite{PC07}}
For a central configuration, the corresponding $\nu$ in the equations $(12)-(17)$ is positive.
\end{lem}

\section{The Proof of Theorem 1.1}
Suppose $q=(q_1,q_2,q_3,q_4)\in (\mathbb{R}^2)^4$ is a planar central configuration as in the hypothesis of Theorem 1.1.

Our goal is to prove $$\triangle_4=-\triangle_3.$$

The way of proving Theorem 1.1 is by the contradiction argument. We assume that $$\triangle_4\neq-\triangle_3.$$

\begin{lem}
 Under the hypothesis of Theorem 1.1, the following inequality holds :
\begin{equation}\triangle_3 <\triangle_1 < 0 < \triangle_2 < \triangle_4. \end{equation}
\end{lem}
\begin{proof}According to the paper \cite{A08}, we know that $|\triangle_{1}|<|\triangle_{3}|$ and $|\triangle_{2}|<|\triangle_{4}|$.
 Because $\triangle_{1}$, $\triangle_{3}$ are negative and $\triangle_{2}$, $\triangle_{4}$ positive, we get $\triangle_3 <\triangle_1 < 0 < \triangle_2 < \triangle_4$.
\end{proof}

\begin{lem} If $\triangle_4\neq -\triangle_3$, then $$b\neq e.$$\end{lem}

\begin{proof}
To prove the Lemma, we need to consider two possible cases:

\textbf{Case 1}. $\triangle_3+\triangle_4 > 0$.

In this case, by the equation (6), we have
\begin{align*}\triangle_1\triangle_3-\triangle_2\triangle_4~&=\triangle_1\triangle_3+(\triangle_1+\triangle_3+\triangle_4)\triangle_4 \\
~&=\triangle_1\triangle_3+\triangle_1\triangle_4+\triangle_3\triangle_4+\triangle_4^2 \\ ~&=(\triangle_1+\triangle_4)(\triangle_3+\triangle_4).
\end{align*}
We claim $$\triangle_1+\triangle_4>0.$$
In fact, if $\triangle_1+\triangle_4\leqq0,$ then $$\triangle_4\leqq-\triangle_1.$$
From the Lemma 3.1, we get $$0<\triangle_2<\triangle_4\leqq-\triangle_1<-\triangle_3.$$
It means $$\triangle_2+\triangle_3 < 0,~~~~ \triangle_1+\triangle_4\leqq0,$$
thus $$\triangle_1+\triangle_2+\triangle_3+\triangle_4 < 0,$$
which contradicts with the equation (6).

So, we have $$\triangle_1\triangle_3-\triangle_2\triangle_4=(\triangle_1+\triangle_4)(\triangle_3+\triangle_4)>0.$$

Thus we get $$\triangle_1\triangle_3 > \triangle_2\triangle_4.$$
By $\alpha > 0$, we then obtain $$\frac{\triangle_1\triangle_3}{\alpha} > \frac{\triangle_2\triangle_4}{\alpha}.$$
Since $\nu > 0$, thus by the equations (13), (16) and the monotonicity of $\varphi'(s)$, we have $$b > e.$$

\textbf{Case 2}.  $\triangle_3+\triangle_4 < 0$.

In this case, the equation (6) implies that  $$\triangle_1\triangle_3-\triangle_2\triangle_4=(\triangle_1+\triangle_4)(\triangle_3+\triangle_4).$$
Since $\triangle_1+\triangle_4>0$, we obtain $$(\triangle_1+\triangle_4)(\triangle_3+\triangle_4)<0.$$
Thus we get $$\triangle_1\triangle_3 < \triangle_2\triangle_4.$$
By $\alpha > 0$, we then obtain $$\frac{\triangle_1\triangle_3}{\alpha} < \frac{\triangle_2\triangle_4}{\alpha}.$$
Since $\nu > 0$, thus similar to Case 1, we obtain $$b < e.$$

Therefore, in both cases, Lemma 3.2 holds.
\end{proof}

From Lemma 3.2, we get
\begin{lem} If $r_{13} = r_{24}$, then $$\triangle_4= -\triangle_3.$$\end{lem}

\begin{proof} We assume $$\triangle_4\neq -\triangle_3.$$
Using Lemma 3.2, we obtain $$b\neq e.$$
It means $r_{13} \neq r_{24}$ and get a contradiction. So $\triangle_3+\triangle_4 \neq 0$ is impossible.
Finally, we get $\triangle_4=-\triangle_3$.

\end{proof}

\begin{lem}If $\triangle_4=-\triangle_3$, then the quadrilateral $q$ is a trapezoid. \end{lem}
\begin{center}
\setlength{\unitlength}{7cm}
\scalebox{0.6}{
\begin{picture}(1.5,1)
\put(0,0){\line(1,0){1.444}}
\put(0,0){\line(1,1){1.04}}
\put(0,0){\line(2,5){.414}}
\put(0.414,1.035){\line(1,0){0.63}}
\put(0.414,1.035){\line(1,-1){1.04}}
\put(1.044,1.035){\line(2,-5){0.414}}
\put(-0.1,0){$q_1$}
\put(1.47,0){$q_2$}
\put(1.044,1.035){$q_3$}
\put(0.314,1.035){$q_4$}

\end{picture}
}

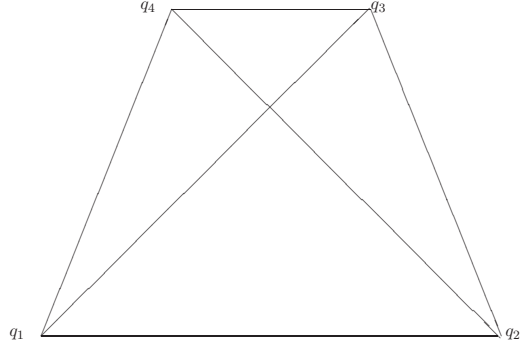
{\bf{Figure.2}}
\end{center}

\begin{proof} Since $\triangle_4=-\triangle_3$, the areas of triangle $\triangle q_{1}q_{2}q_{3}$ and $\triangle q_{1}q_{2}q_{4}$ are equal. We also notice that the two triangles have a common edge $r_{12}$, thus their heights are equal. In other words, the quadrilateral $q$ is a trapezoid.
\end{proof}
Now we can complete the proof of Theorem 1.1:
\begin{proof}[\textbf{Proof of Theorem 1.1}]Under the hypothesis of Theorem 1.1, we have $\triangle_4=-\triangle_3$ from Lemma 3.2 and Lemma 3.3.

 By the equation (6), we get $\triangle_2=-\triangle_1$. It is clear that the configuration is a trapezoid.

By the equations (12)-(17) for the central configurations we have
\begin{equation} \varphi'(r^2_{14})= \frac{\nu}{\alpha} \triangle_1 \triangle_4 + \mu= \frac{\nu}{\alpha} \triangle_2 \triangle_3 + \mu=\varphi'(r^2_{23}), \end{equation}
Since $\varphi'(s)$ is an increasing function of $s$, we obtain $$r_{14}=r_{23}.$$
Therefore, the configuration is an isosceles trapezoid.
\end{proof}

\section{The Proof of Theorem 1.2}

The way to prove Theorem 1.2 is similar to Theorem 1.1.
Under the hypothesis of Theorem 1.2, Lemma 3.1 holds.
Similar to Lemma 3.2, we have the following Lemma:

\begin{lem} If $\triangle_4\neq -\triangle_3$, then $$d\neq c.$$\end{lem}

\begin{proof}
To prove the lemma, we need to consider two possible cases:

\textbf{Case 1}. $\triangle_3+\triangle_4 > 0$.

In this case, by the equation (5), we have $$\triangle_2+\triangle_4 > 0.$$
Using the equation (6), we obtain \begin{align*}\triangle_2\triangle_3-\triangle_1\triangle_4~&=\triangle_2\triangle_3+(\triangle_2+\triangle_3+\triangle_4)\triangle_4 \\
~&=\triangle_2\triangle_3+\triangle_2\triangle_4+\triangle_3\triangle_4+\triangle_4^2 \\ ~&=(\triangle_2+\triangle_4)(\triangle_3+\triangle_4)>0.
\end{align*}

Thus we get $$\triangle_2\triangle_3 > \triangle_1\triangle_4.$$
By $\alpha > 0$, we obtain $$\frac{\triangle_2\triangle_3}{\alpha} > \frac{\triangle_1\triangle_4}{\alpha}.$$
Since $\nu > 0$, thus by equations (14), (15) and the monotonicity of $\varphi'(s)$, we have $$d > c.$$

\textbf{Case 2}.  $\triangle_3+\triangle_4 < 0$.

In this case, the equation (5) implies that $$\triangle_2+\triangle_4 > 0.$$
By the equation (6), we have $$\triangle_2\triangle_3-\triangle_1\triangle_4=(\triangle_2+\triangle_4)(\triangle_3+\triangle_4)<0.$$
Thus we get $$\triangle_2\triangle_3 < \triangle_1\triangle_4.$$
By $\alpha > 0$, we then obtain $$\frac{\triangle_2\triangle_3}{\alpha} < \frac{\triangle_1\triangle_4}{\alpha}.$$
Since $\nu > 0$, thus similar to Case 1, we obtain $$d < c.$$

Therefore, in both cases, Lemma 4.1 holds.
\end{proof}

Using the above Lemma, we get
\begin{lem} If $r_{23} = r_{14}$, then $$\triangle_4= -\triangle_3.$$\end{lem}
\begin{proof} We assume $$\triangle_4\neq -\triangle_3.$$
Using Lemma 4.1, we obtain $$d\neq c.$$
It means $r_{23} \neq r_{14}$ and get a contradiction. So $\triangle_3+\triangle_4 \neq 0$ is impossible.
Finally, we get $\triangle_4=-\triangle_3$.

\end{proof}

\begin{proof}[\textbf{Proof of Theorem 1.2}]Under the hypothesis of Theorem 1.2, we have $\triangle_4=-\triangle_3$ from Lemma 4.1 and Lemma 4.2.

 By the equation (6), we get $\triangle_2=-\triangle_1$. So the configuration ${q}$ is a trapezoid.
 Since $$r_{14}=r_{23},$$
 we obtain that the configuration ${q}$ is an isosceles trapezoid.
\end{proof}

\
\section*{Acknowledgement}
The authors sincerely thank the supports of NSF of China. The authors would like to thank Professor E.Perez-Chavela for his interesting on this paper and some discussions.


\section*{References}
\begin{thebibliography}{10}

\bibitem{A95}  A.Albouy, Symetrie des configurations centrales de quarte corps,
\newblock{\em C.R.Acad.Sci.Paris} {320} (1995) 217-220.

\bibitem{A96}  A.Albouy, The symmetric central configurations of four equal masses,\newblock{\em Contemp. Math } 198 (1996) 131-135.

\bibitem{A08}  A.Albouy, Y.Fu and S.Sun, Symmetry of planar four-body convex central configurations,
\newblock{\em Pro.R.Soc.A} 464 (2008) 1355-1365.


\bibitem{A03} A.Albouy, On a paper of Moeckel on central configurations, \newblock {\em Regul. Chaoic Dyn.} 8 (2003) 133-142.

\bibitem{Cl14} M.Corbera and J.Lliber, Central configuration of the 4-body problem with masses $m_1=m_2>m_3=m_4>0$ and $m$ samll,
\newblock{\em Applied Mathematics and Computation} 246 (2014) 121-147.

\bibitem{CR12} M.Cors and E.Roberts, Four-body co-circular central configurations,
\newblock{\em Nonlinearity} 25 (2012) 343-370.

\bibitem{D19} O.Dziobek, \"{U}ber einen merkw\"{u}rdigen Fall des Viek\"{o}rperproblems,
\newblock{\em Astron.Nachr} 152 (1900) 32-46.

\bibitem{L02} Y.Long and S.Sun, Four-body central configurations with some equal masses,
\newblock{\em Arch Rational Mech.Anal} 162 (2002) 25-44.

\bibitem{M01} R.Moeckel, Generic finitness for Dzibek configurations,
\newblock{\em Trans.Am.Math.Soc} 353 (2001) 4673-4686.

\bibitem{PC07} E.Perez-Chavela and M.Santoprete, Convex four-body central configuratios with some equal masses,
\newblock{\em Arch Rational Mech.Anal} 185 (2007) 481-494.

\bibitem{DS02} D.Schmidt, Central configurations and relative equilibria for the N-body problems, Classical and Celestial Mechanics,
\newblock{\em Princeton, NJ: Princeton University Press} pp 1-33.

\bibitem{DS80} D.G.Saari, On the role and properties of n-body central configurations,
\newblock{\em Celestial Mechanics} 18 (1980) 9-20.

\bibitem{W41} A.Wintner, The analytical foundations of Celestial Mechanics,
\newblock{\em Princeton Math, Series 5, Princeton University Press, Princeton NJ} (1941).


\end{thebibliography}
\end{document}